\documentclass[10pt]{IEEEtran}

\makeatletter
\def\ps@headings{%
\def\@oddhead{\mbox{}\scriptsize\rightmark \hfil \thepage}%
\def\@evenhead{\scriptsize\thepage \hfil \leftmark\mbox{}}%
\def\@oddfoot{}%
\def\@evenfoot{}}
\makeatother
\usepackage{graphicx,amsmath,amssymb,amsfonts,bm,setspace,amsthm}
\usepackage[noadjust]{cite}
\usepackage{braket} 		

\allowdisplaybreaks

\newtheorem{cor}{Corollary}

\newtheorem{lem}{Lemma}

\newtheorem{thm}{Theorem}

\newcommand{\expect}[1]{\mathbb{E}\left[#1\right]}

\newcommand{\prob}[1]{\Pr\left[#1\right]}

\newcommand{\Z}{\mathbb{Z}}
\newcommand{\N}{\mathbb{N}}

\newcommand{\conv}[1]{\operatorname{conv}\left(#1\right)}

\newcommand{\oy}{\overline{y}}
\newcommand{\by}{\bm{y}}

\newcommand{\boy}{\overline{\bm{y}}}
\newcommand{\bor}{\overline{\bm{r}}}

\newcommand{\oor}{\overline{r}}
\newcommand{\bomega}{\bm{\omega}}
\newcommand{\Lambdaint}{\Lambda_{\text{int}}}

\newcommand{\bphi}{\bm{\phi}}
\newcommand{\RRphi}{\mathsf{RR}(\bphi)}
\newcommand{\randrr}{\mathsf{RandRR}}
\newcommand{\sP}{\mathsf{P}}
\newcommand{\blambda}{\bm{\lambda}}

\newcommand{\bmu}{\bm{\mu}}
\newcommand{\bbeta}{\bm{\eta}}
\newcommand{\cvec}[1]{\begin{bmatrix} #1 \end{bmatrix}}
\newcommand{\bv}{\bm{v}}

\newcommand{\sON}{\mathsf{ON}}
\newcommand{\sOFF}{\mathsf{OFF}}

\newcommand{\bQ}{\bm{Q}}
\newcommand{\br}{\bm{r}}
\newcommand{\qrrnum}{\mathsf{QRRNUM}}
\newcommand{\qrr}{\mathsf{QRR}}
\newcommand{\gmax}{G_{\text{max}}}

\newcommand{\tmax}{T_{\text{max}}}

\begin{document}

\title{Network Utility Maximization over Partially Observable Markovian Channels}



\author{\large{Chih-ping~Li,~\IEEEmembership{Student Member,~IEEE} and Michael~J.~Neely,~\IEEEmembership{Senior Member,~IEEE}}%
\thanks{Chih-ping Li (web: http://www-scf.usc.edu/$\sim$chihpinl) and Michael J. Neely (web:  http://www-rcf.usc.edu/$\sim$mjneely) are with the Department of Electrical Engineering, University of Southern California, Los Angeles, CA 90089, USA.}    %
\thanks{This material is supported in part  by one or more of the following: the DARPA IT-MANET program grant W911NF-07-0028, the NSF Career grant CCF-0747525, and continuing through participation in the Network Science Collaborative Technology Alliance sponsored by the U.S. Army Research Laboratory.}%
}


\maketitle

\begin{abstract}
We consider a utility maximization problem over partially observable Markov $\sON$/$\sOFF$ channels. In this network instantaneous channel states are never known, and at most one user is selected for service in every slot according to the partial channel information provided by past observations. Solving the utility maximization problem directly is difficult because it involves solving partially observable Markov decision processes. Instead, we construct an approximate solution by optimizing the network utility only over a good constrained network capacity region rendered by stationary policies. Using a novel frame-based Lyapunov drift argument, we design a policy of admission control and user selection that stabilizes the network with utility that can be made arbitrarily close to the optimal in the constrained region. Equivalently, we are dealing with a high-dimensional restless bandit problem with a general functional objective over Markov $\sON$/$\sOFF$ restless bandits. Thus the network control algorithm developed in this paper serves as a new approximation methodology to attack such complex restless bandit problems.

\end{abstract}

\section{Introduction} \label{sec:intro}
 
\begin{figure}[htbp]
\centering
\includegraphics[width=3in]{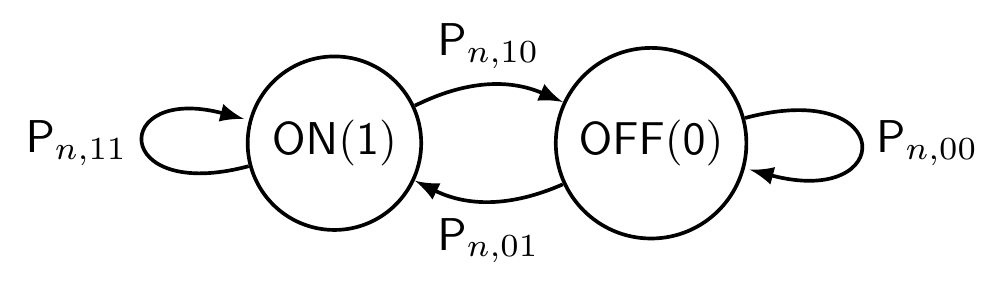}
%
%
\caption{The Markov $\sON$/$\sOFF$ chain for channel $n\in\{1, 2, \ldots, N\}$.}
\label{fig:102}
\end{figure}

This paper studies a multi-user wireless scheduling problem over partially observable environments. We consider a wireless uplink  system serving $N$ users via $N$ independent Markov $\sON$/$\sOFF$ channels (see Fig.~\ref{fig:102}). 
Suppose time is slotted with normalized slots $t\in\Z^{+}$. Channel states are fixed in every slot, and can only change at slot boundaries. In every slot, the channel states are unknown, and at most one user is selected for transmission. The chosen user can successfully deliver a packet  if the channel is $\sON$, and zero otherwise. Since channels are $\sON$/$\sOFF$, the state of the used channel is uncovered by an error-free ACK/NACK feedback at the end of the slot (failing to receive an ACK is regarded as a NACK). The states of each Markovian channel are correlated over time, and thus the revealed channel condition from ACK/NACK feedback provides partial information of future states, which can be used to improve user selection decisions and network performance. Our goal is to design a network control policy that maximizes a general network utility metric which is a  function of the achieved throughput vector. Specifically, let $y_{n}(t)$ be the amount of user-$n$ data served in slot $t$, and define the throughput $\oy_{n}$ for user $n$ as
$
\oy_{n} \triangleq \lim_{t\to\infty} \frac{1}{t} \sum_{\tau=0}^{t-1} \expect{y_{n}(\tau)}
$. 
Let $\Lambda$ be the \emph{network capacity region} of the wireless uplink, defined as the closure of the set of all achievable throughput vectors $\boy\triangleq (\oy_{n})_{n=1}^{N}$. Then we seek to solve the following utility maximization problem:
\begin{align}
\text{maximize:} \quad & g(\boy) \label{eq:101} \\
\text{subject to:} \quad & \boy \in\Lambda \label{eq:102}
\end{align}
where in the above we denote by $g(\cdot)$ a generic utility function that is concave, continuous, nonnegative, and nondecreasing.

The problem~\eqref{eq:101}-\eqref{eq:102} is very important to explore because it has many applications in various fields. In multi-user wireless scheduling, optimizing  network utility over stochastic networks is first solved in~\cite{NML08}, under the assumption that channel states are i.i.d. over slots and are known perfectly and instantly. The problem~\eqref{eq:101}-\eqref{eq:102} we consider here generalizes the network utility maximization framework in~\cite{NML08} to networks with limiting channel probing capability (see \cite{LaN10-DCA, CPA09conf} and references therein) and delayed/uncertain channel state information (see \cite{PET07conf,YaS08conf,YaS09conf} and references therein), in which  we shall take advantage of channel memory~\cite{ZRM96conf} to improve network performance. 
In sequential decision making, \eqref{eq:101}-\eqref{eq:102} also captures an important class of restless bandit problems~\cite{Whi88} in which each Markovian channel represents a two-state restless bandit, and packets served over a channel are rewards from playing the bandit. This class of Markov $\sON$/$\sOFF$ restless bandit problems  has modern applications in opportunistic spectrum access in cognitive radio networks~\cite{ZaS07_1,ZaS07_2} and target tracking of unmanned aerospace vehicles~\cite{NDF08conf}.


Solving the maximization problem~\eqref{eq:101}-\eqref{eq:102} is difficult because $\Lambda$ is unknown. In principle, we may compute $\Lambda$ by locating its boundary points. However, they are solutions to $N$-dimensional Markov decision processes with information state vectors $\bomega(t) \triangleq (\omega_{n}(t))_{n=1}^{N}$, where $\omega_{n}(t)$ is the conditional probability that  channel $n$ is $\sON$ in slot $t$ given the channel observation history. Namely, let $s_{n}(t)$ denote the state of channel $n$ in slot $t$. Then
\begin{equation} \label{eq:112}
\omega_{n}(t) \triangleq \prob{s_{n}(t)=\sON \mid \text{channel observation history}}.
\end{equation}
We will show later $\omega_{n}(t)$ takes values in a countably infinite set. Thus computing $\Lambda$ and solving~\eqref{eq:101}-\eqref{eq:102} seem to be infeasible.

Instead of solving~\eqref{eq:101}-\eqref{eq:102}, in this paper we adopt an achievable region approach to construct approximate solutions to~\eqref{eq:101}-\eqref{eq:102}. The key idea is two-fold. First, we explore the problem structure and construct an achievable throughput region $\Lambdaint \subset \Lambda$ rendered by \emph{good} stationary (possibly randomized) policies. Then we solve the constrained maximization problem:
\begin{align}
\text{maximize:} \quad & g(\boy) \label{eq:104} \\
\text{subject to:} \quad &  \boy \in\Lambdaint \label{eq:105}
\end{align}
as an approximation to~\eqref{eq:101}-\eqref{eq:102}. This approximation is practical because  every throughput vector in $\Lambdaint$ is attainable by simple stationary policies, and achieving feasible points outside $\Lambdaint$ may require solving the much more complicated partially observable Markov decision processes (POMDPs) that relate to the original problem. Thus for the sake of  simplicity and practicality, we shall regard $\Lambdaint$ as our \emph{operational} network capacity region.

Using the rich structure of the Markovian channels, in~\cite{LaN10conf-channelmemory, LaN10arXiv-channelmemory} we  have constructed a good achievable region $\Lambdaint$ rendered by a special class of \emph{randomized round robin}  policies. It is important to note that we will maximize $g(\boy)$ only over this class of policies. Since every point in $\Lambdaint$ can be achieved by one such policy (which we will show later), equivalently we are solving~\eqref{eq:104}-\eqref{eq:105}. We remark that solving~\eqref{eq:104}-\eqref{eq:105} is decoupled from the construction of $\Lambdaint$. We will show in this paper that~\eqref{eq:104}-\eqref{eq:105} can be solved. Therefore, the overall optimality of this achievable region approach depends on the proximity of the inner bound $\Lambdaint$ to the full capacity  region $\Lambda$.

The main contribution of this paper is that, using the Lyapunov optimization theory originally developed in~\cite{TaE92, TaE93} and later generalized by~\cite{Nee03thesis, NML08} for optimal stochastic control over wireless networks (see~\cite{GNT06} for an introduction), we can solve~\eqref{eq:104}-\eqref{eq:105} and develop optimal greedy algorithms. Specifically, using a novel Lyapunov drift argument, we construct a frame-based, queue-dependent network control algorithm of service allocation and admission control.\footnote{Admission control is used to facilitate the solution to the problem~\eqref{eq:104}-\eqref{eq:105}.} At the beginning of each frame, the admission controller decides how much new data to admit by solving a simple convex program.\footnote{The admission control decision decouples into $N$ separable one-dimensional problems that are easily solved in real time in the case when $g(\boy)$ is a sum of one-dimensional utility functions for each user.}  The service allocation decision selects a randomized round robin policy by maximizing an \emph{average MaxWeight} metric, and runs the policy for one round in the frame. We will show that this joint policy stabilizes the network and yields the achieved network utility $g(\boy)$ satisfying
\begin{equation} \label{eq:710}
g(\boy) \geq g(\boy^{*}) - \frac{B}{V_{g}},
\end{equation}
where $g(\boy^{*})$ is the optimal objective of~\eqref{eq:104}-\eqref{eq:105}, $B>0$ is a finite constant, $V_{g}$ is a predefined positive control parameter, and we temporarily assume that all limits exist. By choosing $V_{g}$ sufficiently large, we can approach the optimal  utility $g(\boy^{*})$ arbitrarily well in~\eqref{eq:710}, and thus solve~\eqref{eq:104}-\eqref{eq:105}.


Restless bandit problems with Markov $\sON$/$\sOFF$ bandits  have been studied in~\cite{LaZ08arXiv,Nin08conf,GMS09arXiv}, in which \emph{index policies}~\cite{Whi88,Nin07} are developed to maximize long-term average/discounted rewards. In this paper we extend this class of problems to having a general functional objective that needs to be maximized. This new problem is difficult to solve using existing approaches such as Whittle's index~\cite{Whi88} or Markov decision theory~\cite{Yu06thesis}, because they are typically  limited to deal with problems with very simple objectives. The achievable region approach we develop in this paper solves (approximately) this extended problem, and thus could be viewed as a new approximation methodology  to analyze similar complex restless bandit problems.

In the next section we introduce the detailed network model. Section~\ref{sec:inn} summarizes the construction of the inner bound $\Lambdaint$ in~\cite{LaN10conf-channelmemory,LaN10arXiv-channelmemory}. Our dynamic control algorithm is developed in Section~\ref{sec:num}, and the performance analysis is given in Section~\ref{sec:analysis}.

\section{Detailed Network Model} \label{sec:admit}
In addition to the basic network model given in Section~\ref{sec:intro},  we suppose every channel $n\in\{1, \ldots, N\}$ evolves according to the transition probability matrix
\[
\bm{\sP}_{n}
=
\begin{bmatrix}
\sP_{n, 00} & \sP_{n, 01} \\
\sP_{n, 10} & \sP_{n, 11}
\end{bmatrix},
\]
where state $\sON$ is represented by $1$ and $\sOFF$ by $0$, and $\sP_{n, ij}$ denotes the transition probability from state $i$ to $j$.  We suppose every channel is \emph{positively correlated over time}, so that an $\sON$ state is likely to be followed by another $\sON$ state. An equivalent mathematical definition is $x_{n} \triangleq \sP_{n, 01} + \sP_{n, 10} <1$ for all $n$. Let $\bm{\sP}_{n}$ be known by both the network and user $n$.

We suppose every user has a data source of unlimited packets.  In every slot,  user $n\in\{1, \ldots, N\}$ admits $r_{n}(t)\in[0,1]$ packets from the source into a queue $Q_{n}(t)$ of infinite capacity. For simplicity, we assume $r_{n}(t)$ takes real values in $[0, 1]$ for all $n$.\footnote{We can accommodate the integer-value assumption of $r_{n}(t)$ by introducing \emph{auxiliary queues}; see~\cite{NML08} for an example.}  Define $\br(t) \triangleq (r_{n}(t))_{n=1}^{N}$. At the beginning of every slot, the network chooses and sends to the users one feasible admitted data vector $\br(t)$ according to some admission policy. We let  $Q_{n}(t)$ and $\mu_{n}(t)\in\{0, 1\}$ denote the queue backlog and the service rate of user $n$ in slot $t$. Assume $Q_{n}(0) = 0$ for all $n$. Then the queueing process $\{Q_{n}(t)\}$ evolves as
\begin{equation} \label{eq:202}
Q_{n}(t+1) = \max[Q_{n}(t) - \mu_{n}(t), 0] + r_{n}(t).
\end{equation}
The network keeps track of the backlog vector $\bQ(t)\triangleq (Q_{n}(t))_{n=1}^{N}$ in every slot. We say queue $Q_{n}(t)$ is (strongly) stable if
\[
\limsup_{t\to\infty} \frac{1}{t} \sum_{\tau=0}^{t-1} \expect{Q_{n}(\tau)} < \infty,
\]
and the network is  stable if all queues in the network are stable. Clearly a sufficient condition for stability is:
\begin{equation} \label{eq:706}
\limsup_{t\to\infty} \frac{1}{t} \sum_{\tau=0}^{t-1} \sum_{n=1}^{N} \expect{Q_{n}(\tau)} < \infty.
\end{equation}
Our goal is to design a policy that admits the right amount of packets into the network and serves them properly, so that the network is stable with utility that can be made arbitrarily close to the optimal solution to~\eqref{eq:104}-\eqref{eq:105}.

\section{A Performance Inner Bound} \label{sec:inn}
In this section we summarize the results in~\cite{LaN10conf-channelmemory,LaN10arXiv-channelmemory} on constructing an achievable region $\Lambdaint$ using randomized round robin policies. 
See~\cite{LaN10arXiv-channelmemory} for detailed proofs.

\subsection{Sufficient statistic}

As discussed in~\cite[Chapter $5.4$]{Ber05book}, the information state vector $\bomega(t)$ defined in~\eqref{eq:112} is a \emph{sufficient statistic} of the network, meaning that it suffices to make optimal decisions based only on $\bomega(t)$ in every slot.

For channel $n\in\{1, \ldots, N\}$, we denote by $\sP_{n,ij}^{(k)}$ the $k$-step transition probability from state $i$ to $j$, and $\pi_{n, \sON}$ its stationary probability of state $\sON$. Since channels are positively correlated, we can show that $\sP_{n, 01}^{(k)}$ is nondecreasing and $\sP_{n, 11}^{(k)}$ is nonincreasing in $k$, and
$
\pi_{n, \sON} = \lim_{k\to\infty} \sP_{n, 01}^{(k)} = \lim_{k\to\infty} \sP_{n, 11}^{(k)}
$.
For channel $n$, conditioning on the outcome of the last observation and when it was taken, it is easy to see that $\omega_{n}(t)$ takes values in the countably infinite set $\mathcal{W}_{n} \triangleq \{\sP_{n, 01}^{(k)}, \sP_{n, 11}^{(k)} : k\in\N\} \cup \{\pi_{n, \sON}\}$. Let $n(t)$ be the channel observed in slot $t$ via ACK/NACK feedback. The evolution of $\omega_{n}(t)$ for each $n$ then follows:
\begin{equation} \label{eq:113}
\omega_{n}(t+1) \!=\!
\begin{cases}
\sP_{n, 01},\quad \text{if $n=n(t)$, $s_{n}(t)=\sOFF$} \\
\sP_{n, 11},\quad \text{if $n=n(t)$, $s_{n}(t)=\sON$} \\
\omega_{n}(t) \sP_{n, 11} + (1-\omega_{n}(t))\sP_{n, 01},\;\text{if $n\neq n(t)$}.
\end{cases}
\end{equation}



\subsection{Randomized round robin} \label{sec:randrr}

Let $\Phi$ denote the set of all $N$-dimensional binary vectors excluding the zero vector $\bm{0}$. Every vector $\bphi \triangleq (\phi_{n})_{n=1}^{N} \in \Phi$ stands for a collection of \emph{active} channels, where we say channel $n$ is active in $\bphi$ if $\phi_{n}=1$.  Let $M(\bphi)$ denote the number of $1$'s (or active channels) in $\bphi$.

Consider the following \emph{dynamic round robin} policy $\RRphi$ that serves active channels in $\bphi$ possibly with different order in different rounds. This is the building block of the randomized round robin policies that we will introduce shortly.

\underline{\bf Dynamic Round Robin Policy $\RRphi$: }
\begin{enumerate}
\item In each round, suppose an ordering of active channels in $\bphi$ is given. \label{item:201}

\item When switching to active channel $n$, with probability $\sP_{n,01}^{(M(\bphi))} / \omega_{n}(t)$ keep transmitting packets  over channel $n$ until a NACK is received, and then switch to the next active channel.  With probability $1- \sP_{n,01}^{(M(\bphi))} / \omega_{n}(t)$, transmit a dummy packet with no information content for one slot (used for channel sensing) and then switch to the next active channel.

\item Update $\bomega(t)$ according to~\eqref{eq:113} in every slot.
\end{enumerate}

It is shown in~\cite{ALJ09} that, when channels have the same transition probability matrix, serving all channels by a greedy round robin policy maximizes the sum throughput of the network. Thus we shall get a good achievable throughput region $\Lambdaint$ by randomly mixing round robin policies, each of which serves a different subset of channels. 

Consider the following randomized round robin that mixes $\RRphi$ policies for different $\bphi$:

\underline{\bf Randomized Round Robin Policy $\randrr$}:
\begin{enumerate}
\item Pick $\bphi \in \Phi \cup \{\bm{0}\}$ with probability $\alpha_{\bphi}$, where $\alpha_{\bm{0}} + \sum_{\bphi\in\Phi} \alpha_{\bphi} = 1$. \label{item:4}
\item If $\bphi\in\Phi$ is selected, run $\mathsf{RR}(\bphi)$ for one round with the channel ordering of \emph{least recently used first}. Then go to Step~\ref{item:4}. If $\bphi = \bm{0}$, idle the system for one slot and then go to Step~\ref{item:4}. \label{item:7}
\end{enumerate}

For notational convenience, let $\mathsf{RR}(\bm{0})$ denote the operation of idling the system for one slot. For any $\bphi\in\Phi$, we note that the $\RRphi$ policy is feasible only if $\sP_{n, 01}^{(M(\bphi))} \leq \omega_{n}(t)$ whenever we switch to active channel $n$. This condition is enforced in every $\randrr$ policy by serving active channels in the order of \emph{least recently used first}~\cite[Lemma $6$]{LaN10arXiv-channelmemory}. Consequently, every $\randrr$ is a feasible policy.\footnote{The feasibility of $\randrr$ policies is proved in~\cite{LaN10arXiv-channelmemory} under the special case that there are no idle operations ($\alpha_{\bm{0}}=0$).  Using the monotonicity of  $k$-step transition probabilities $\{\sP_{n, 01}^{(k)}, \sP_{n, 11}^{(k)}\}$, the feasibility can be similarly proved for the generalized $\randrr$ policies considered here.}  We note that the $\randrr$ policies considered here are a superset of those in~\cite{LaN10arXiv-channelmemory}, because here we allow the additional idling operation. This enlarged policy space, however, has the same achievable throughput region as that in~\cite{LaN10arXiv-channelmemory}, because idle operations do not improve throughput. We generalize the $\randrr$ policies here to ensure that every feasible point in $\Lambdaint$ can be achieved by some $\randrr$ policy.
It is also helpful to note that, for any $\bphi\in\Phi$ and a fixed channel ordering, every $\RRphi$ policy is a special case of the randomized round robin $\randrr$ with $\alpha_{\bphi} = 1$ and $0$ otherwise.


\subsection{The achievable region}

Next we summarize the achievable region rendered by randomized round robin  policies.

\begin{thm}[\cite{LaN10conf-channelmemory,LaN10arXiv-channelmemory}] \label{thm:101}
For each vector $\bphi\in\Phi$, define the $N$-dimensional vector $\bbeta^{(\bphi)} \triangleq (\eta^{(\bphi)}_{n})_{n=1}^{N}$ where
\[
\eta^{(\bphi)}_{n} \triangleq
\begin{cases}
\frac{	
	\sP_{n,01} (1-(1-x_{n})^{M(\bphi)}) / (x_{n} \sP_{n,10})
}{
	M(\bphi) + \sum_{n: \phi_{n}=1} \frac{\sP_{n,01} (1-(1-x_{n})^{M(\bphi)}) }{ x_{n} \sP_{n,10}}
}, & \text{if $\phi_{n} = 1$} \\
0, & \text{if $\phi_{n} = 0$}
\end{cases}
\]
and $x_{n} \triangleq \sP_{n, 01} + \sP_{n,10}$. Then the class of $\mathsf{RandRR}$ policies supports all throughput vectors $\blambda$ in the set
\[
\Lambdaint \triangleq \left\{
	\blambda \mid \bm{0}\leq \blambda \leq \bmu,\ \bmu\in\conv{
	\left\{ \bbeta^{\bphi}\right\}_{\bphi\in\Phi}
}
\right\},
\]
where $\conv{A}$ denotes the convex hull of set $A$, and $\leq$ is taken entrywise.
\end{thm}

\begin{cor} \label{cor:101}
When channels have the same transition probability matrix so that $\bm{\sP}_{n} = \bm{\sP}$ for all $n$, we have:
\[
\Lambda_{\text{int}} = \left\{
\blambda \mid
\bm{0}\leq \blambda \leq \bmu,\,\bmu\in\conv{ \left\{\frac{c_{M(\bphi)}}{M(\bphi)} \bphi \right\}_{\bphi\in\Phi}}
\right\},
\]
where
\begin{equation} \label{eq:107}
c_{M(\bphi)} \triangleq \frac{
	\sP_{01}(1-(1-x)^{M(\bphi)})
}{
	x\,\sP_{10} + \sP_{01}(1-(1-x)^{M(\bphi)})
}, \quad x = \sP_{01} + \sP_{10},
\end{equation}
and we have dropped the subscript $n$ due to channel symmetry.
\end{cor}

The closeness of the inner bound $\Lambdaint$ and the full capacity region $\Lambda$ is quantified in~\cite{LaN10arXiv-channelmemory} in the special case that channels have the same transition probability matrix. For any feasible direction $\bv$, it can be shown that as $\bv$ becomes more symmetric, or forms a smaller angle with the $45$-degree line, the loss of the sum throughput of the inner boundary point  in direction $\bv$ decreases to zero geometrically fast, provided that the network serves a large number of users.

Next, that $\randrr$ policies considered in this paper are random mixings of those in~\cite{LaN10arXiv-channelmemory} and idle operations leads to the next corollary.

\begin{cor} \label{cor:301}
Every throughput vector in $\Lambdaint$ can be achieved by some $\randrr$ policy.
\end{cor}

\subsection{A two-user example}
\begin{figure}[htbp]
\centering
\includegraphics[width=2.6in]{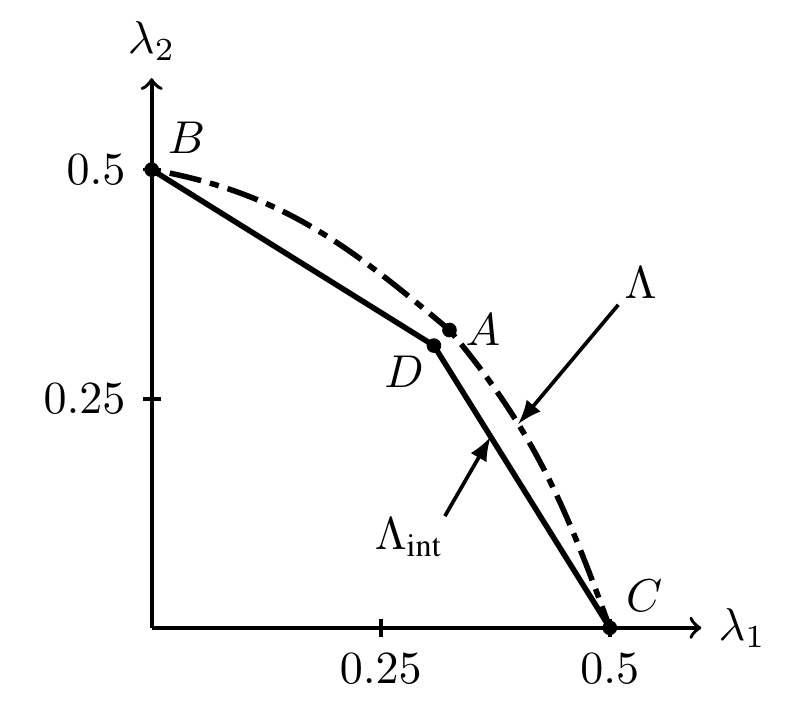}

\caption{The closeness of $\Lambdaint$ and $\Lambda$.}
\label{fig:101}
\end{figure}

Consider a two-user system with symmetric channels with $\sP_{01} = \sP_{10} = 0.2$. From Corollary~\ref{cor:101}
\[
\Lambda_{\text{int}} = \Set{
\cvec{\lambda_{1} \\ \lambda_{2}} |
	\begin{gathered}[c]
		0\leq \lambda_{n}\leq \mu_{n}, \text{ for $1\leq n\leq 2$,} \\
		\cvec{\mu_{1}\\ \mu_{2}}\in\conv{
			\left\{
			\cvec{c_{2}/2 \\ c_{2}/2},
			\cvec{c_{1} \\ 0},
			\cvec{0 \\ c_{1}}
			\right\}
		}
	\end{gathered}
}.
\]
where $c_{1}$ and $c_{2}$ are defined in~\eqref{eq:107}. Fig.~\ref{fig:101} shows the closeness of $\Lambdaint$ and $\Lambda$ in this example.
We note that points $B$, $A$, and $C$ maximize the sum throughput of the network in directions $(0,1)$, $(1,1)$, and $(1,0)$, respectively~\cite{ALJ09}. Therefore the boundary of the (unknown)  full capacity region $\Lambda$ is a concave curve connecting these points.

\section{Network Utility Maximization} \label{sec:num}

From Theorem~\ref{thm:101},  the constrained problem~\eqref{eq:104}-\eqref{eq:105} is a well-defined convex program. However, solving~\eqref{eq:104}-\eqref{eq:105} remains difficult because the representation of $\Lambdaint$ via a convex hull of $(2^{N}-1)$ throughput vectors is very complicated. Next we solve~\eqref{eq:104}-\eqref{eq:105} by admission control and service allocation in the network. We will use the Lyapunov optimization theory to construct a dynamic policy that \emph{learns} a near-optimal solution to~\eqref{eq:104}-\eqref{eq:105}, where the closeness to the true optimality is controlled by a positive control parameter $V_{g}$.


\subsection{Constructing Lyapunov drift} \label{sec:drift}

We start with constructing a frame-based Lyapunov \emph{drift-minus-utility} inequality over a frame of size $T$, where $T$ is possibly random but has a finite second moment bounded by a constant $C$ so that $C\geq \expect{T^{2}\mid \bQ(t)}$ for all $t$ and all possible $\bQ(t)$. Define $B\triangleq NC$. The result will shed light on the structure of our desired policy. By iteratively applying~\eqref{eq:202}, it is not hard to show that
\begin{equation} \label{eq:203}
Q_{n}(t+T) \leq \max \left[Q_{n}(t) - \sum_{\tau=0}^{T-1} \mu_{n}(t+\tau) , 0\right] + \sum_{\tau=0}^{T-1} r_{n}(t+\tau)
\end{equation}
 for each $n\in\{1, \ldots, N\}$. We define the \emph{Lyapunov~function}
 \[
 L(\bQ(t)) \triangleq \frac{1}{2} \sum_{n=1}^{N} Q_{n}^{2}(t)
 \]
 and the \emph{$T$-slot Lyaupnov drift}
\[
\Delta_{T}(\bQ(t)) \triangleq \expect{L(\bQ(t+T)) - L(\bQ(t)) \mid \bQ(t)},
\]
where the expectation is over the randomness of the network in the frame, including that of $T$.  By taking the following steps: (1) take square of~\eqref{eq:203} for each $n$; (2) use the inequalities
\begin{gather*}
\max[a-b, 0] \leq a, \quad \forall a\geq 0, \\
(\max[a-b, 0])^{2} \leq (a-b)^{2}, \quad \mu_{n}(t) \leq 1, \quad r_{n}(t) \leq 1,
\end{gather*}
to simplify terms; (3) sum all resulting inequalities; (4) take conditional expectation on $\bQ(t)$, we can show
\begin{equation} \label{eq:204}
\begin{split}
& \Delta_{T}(\bQ(t)) \leq B \\
& - \expect{ \sum_{n=1}^{N} Q_{n}(t) \left[ \sum_{\tau=0}^{T-1} \mu_{n}(t+\tau) - r_{n}(t+\tau)\right] \mid \bQ(t)}.
\end{split}
\end{equation}
By subtracting from both sides of~\eqref{eq:204} the weighted sum utility
\[
V_{g} \expect{ \sum_{\tau=0}^{T-1} g(\br(t+\tau)) \mid \bQ(t)},
\]
where $V_{g}>0$ is a predefined control parameter, we get
\begin{align}
& \Delta_{T}(\bQ(t)) - V_{g} \expect{ \sum_{\tau=0}^{T-1} g(\br(t+\tau)) \mid \bQ(t)} \notag \\
& \leq B - \sum_{n=1}^{N} Q_{n}(t) \expect{\sum_{\tau=0}^{T-1} \mu_{n}(t+\tau) \mid \bQ(t)} \notag \\
& - \expect{ \sum_{\tau=0}^{T-1} \left[V_{g}\,g(\br(t+\tau)) - \sum_{n=1}^{N} Q_{n}(t) r_{n}(t+\tau) \right]\mid \bQ(t)}. \label{eq:205}
\end{align}
The above inequality gives an upper bound on the \emph{drift-minus-utility} expression at the left side of~\eqref{eq:205}, and holds for any scheduling policy over a frame of any size $T$.

\subsection{Network control policy}

Let $f(\bQ(t))$ and $g(\bQ(t))$ denote the second-to-last and the last term of~\eqref{eq:205}:
\begin{align*}
f(\bQ(t)) &\triangleq \sum_{n=1}^{N} Q_{n}(t) \expect{\sum_{\tau=0}^{T-1} \mu_{n}(t+\tau) \mid \bQ(t)} \\
g(\bQ(t)) &\triangleq \mathbb{E}\Bigg[ \sum_{\tau=0}^{T-1} \Bigg[V_{g}\,g(\br(t+\tau)) \\
	&\qquad\qquad\qquad\quad - \sum_{n=1}^{N} Q_{n}(t) r_{n}(t+\tau) \Bigg]\mid \bQ(t)\Bigg],
\end{align*}
and~\eqref{eq:205} is equivalent to
\begin{equation} \label{eq:309}
\begin{split}
&\Delta_{T}(\bQ(t)) - V_{g} \expect{ \sum_{\tau=0}^{T-1} g(\br(t+\tau)) \mid \bQ(t)} \\
&\quad \leq B - f(\bQ(t)) - g(\bQ(t)).
\end{split}
\end{equation}
After observing the current backlog vector $\bQ(t)$, we seek to maximize over all feasible policies the average
\begin{equation} \label{eq:401}
\frac{f(\bQ(t))+g(\bQ(t))}{\expect{T \mid \bQ(t)}}
\end{equation}
over a frame of size $T$. Every feasible policy here consists of: (1) an admission policy that admits $r_{n}(t+\tau)$ packets to user $n$ in every slot of the frame, and (2) a randomized round robin $\randrr$ policy introduced in Section~\ref{sec:randrr} that serves a set of active users and decides the service rates $\mu_{n}(t+\tau)$ in the frame. The random frame size $T$ in~\eqref{eq:401} is  the length of one transmission round under the candidate $\randrr$ policy, and its distribution depends on the backlog vector $\bQ(t)$ via the queue-dependent choice of $\randrr$. We will show later that the novel performance metric~\eqref{eq:401} helps to achieve near-optimal network utility.

We simplify the procedure of maximizing~\eqref{eq:401} in the following, and the result is our network control algorithm. In $g(\bQ(t))$, we observe that the optimal choices of the admitted data vectors $\br(t+\tau)$ are independent of both the frame size $T$ and the rate allocations $\mu_{n}(t+\tau)$ in $f(\bQ(t))$. Thus, $\br(t+\tau)$ can be optimized separately. Specifically, the optimal values of $\br(t+\tau)$ shall be the same for all $\tau\in\{0, \ldots, T-1\}$ and are the solution to
\begin{align}
\text{maximize:} & \quad  V_{g}\,g(\br(t))  - \sum_{n=1}^{N} Q_{n}(t) r_{n}(t) \label{eq:209}\\
\text{subject to:} & \quad r_{n}(t) \in [0, 1], \quad \forall n\in\{1, \ldots, N\} \label{eq:210}
\end{align}
which only depends on the backlog vector $\bQ(t)$ at the beginning of the current frame and the predefined control parameter $V_{g}$. We note that if $g(\cdot)$ is a sum of individual utilities so that $g(\br(t)) = \sum_{n=1}^{N} g_{n}(r_{n}(t))$, \eqref{eq:209}-\eqref{eq:210} decouples into $N$ one-dimensional convex programs, each of which maximizes $V_{g}\,g_{n}(r_{n}(t))  -  Q_{n}(t) r_{n}(t)$ over $r_{n}(t) \in [0, 1]$, which can be solved efficiently in real time. Let $h^{*}(\bQ(t))$ be the resulting optimal objective of~\eqref{eq:209}-\eqref{eq:210}. It follows that
\[
g(\bQ(t)) = \expect{T \mid \bQ(t)}  h^{*}(\bQ(t))
\]
and~\eqref{eq:401} is equal to
\begin{equation} \label{eq:402}
\frac{f(\bQ(t))}{\expect{T \mid \bQ(t)}} + h^{*}(\bQ(t)).
\end{equation}
The sum~\eqref{eq:402} indicates that finding the optimal admission policy is independent of finding the optimal randomized round robin policy. It remains to maximize the first term of~\eqref{eq:402} over all $\randrr$ policies. 

Next we evaluate the first term of~\eqref{eq:402} under a fixed $\randrr$ policy with parameters $\{\alpha_{\bphi}\}_{\bphi\in\Phi \cup \{\bm{0}\}}$. Conditioning on the choice of $\bphi$, we get
\[
f(\bQ(t)) = \sum_{\bphi\in\Phi\cup\{\bm{0}\}} \alpha_{\bphi} f(\bQ(t), \RRphi),
\]
where $f(\bQ(t), \RRphi)$ denotes the term $f(\bQ(t))$ evaluated under policy $\RRphi$ (recall that $\RRphi$ is a special case of the $\randrr$ policy). Similarly, by conditioning we can show
\[
\expect{T} = \expect{T \mid \bQ(t)} = \sum_{\bphi\in\Phi\cup\{\bm{0}\}} \alpha_{\bphi} \expect{T_{\RRphi}},\footnote{Given a fixed policy $\randrr$, the frame size $T$ no longer depends on the backlog vector $\bQ(t)$. Therefore $\expect{T} = \expect{T\mid\bQ(t)}$.}
\]
where $T_{\RRphi}$ denotes the duration of one transmission round under the $\RRphi$ policy. It follows that
\begin{equation} \label{eq:403}
\frac{f(\bQ(t))}{\expect{T \mid \bQ(t)}} = \frac{\sum_{\bphi\in\Phi\cup\{\bm{0}\}} \alpha_{\bphi} f(\bQ(t), \RRphi)}{\sum_{\bphi\in\Phi\cup\{\bm{0}\}} \alpha_{\bphi} \expect{T_{\RRphi}}}.
\end{equation}

The next lemma shows there always exists a $\RRphi$ policy maximizing~\eqref{eq:403} over all $\randrr$ policies. Therefore it suffices to focus only on $\RRphi$ policies.

\begin{lem} \label{lem:401}
We index $\RRphi$ policies for all $\bphi\in\Phi\cup\{\bm{0}\}$. For the $\RRphi$ policy with index $k$, define 
\[
f_{k} \triangleq f(\bQ(t), \RRphi), \quad D_{k} \triangleq \expect{T_{\RRphi}}.
\]
Without loss of generality, assume
\[
\frac{f_{1}}{D_{1}} \geq \frac{f_{k}}{D_{k}}, \quad \forall k\in\{2, 3, \ldots, 2^{N}\}.
\]
Then for any probability distribution $\{\alpha_{k}\}_{k\in\{1, \ldots, 2^{N}\}}$ with $\alpha_{k}\geq 0$ and $\sum_{k} \alpha_{k}=1$, we have
\[
\frac{f_{1}}{D_{1}} \geq \frac{\sum_{k=1}^{2^{N}} \alpha_{k} f_{k}}{\sum_{k=1}^{2^{N}} \alpha_{k} D_{k}}.
\]
\end{lem}

\begin{IEEEproof}[Proof of Lemma~\ref{lem:401}]
Fact $1$: Let $\{a_{1}, a_{2}, b_{1}, b_{2}\}$ be four positive numbers, and suppose there is a bound $z$ such that $a_{1}/b_{1} \leq z$ and $a_{2}/b_{2} \leq z$. Then for any probability $\theta$ (where $0\leq \theta\leq 1$), we have:
\begin{equation} \label{eq:404}
\frac{\theta a_{1}+ (1-\theta)a_{2}}{\theta b_{1}+ (1-\theta)b_{2}} \leq z.
\end{equation}

We prove Lemma~\ref{lem:401} by induction and~\eqref{eq:404}. Initially, for any $\alpha_{1}$, $\alpha_{2}\geq 0$, $\alpha_{1}+\alpha_{2}=1$, from $f_{1}/D_{1} \geq f_{2}/D_{2}$ we get
\[
\frac{f_{1}}{D_{1}} \geq \frac{\alpha_{1} f_{1}+\alpha_{2} f_{2}}{\alpha_{1} D_{1}+\alpha_{2} D_{2}}.
\]
For some $K > 2$, assume
\begin{equation} \label{eq:405}
\frac{f_{1}}{D_{1}} \geq \frac{\sum_{k=1}^{K-1} \alpha_{k} f_{k}}{\sum_{k=1}^{K-1} \alpha_{k} D_{k}}
\end{equation}
holds for any probability distribution $\{\alpha_{k}\}_{k=1}^{K-1}$. It follows that, for any probability distribution $\{\alpha_{k}\}_{k=1}^{K}$, we get
\[
\frac{\sum_{k=1}^{K} \alpha_{k} f_{k}}{\sum_{k=1}^{K} \alpha_{k} D_{k}} \!=\! \frac{(1-\alpha_{K})\! \left[\sum_{k=1}^{K-1} \frac{\alpha_{k}}{1-\alpha_{K}} f_{k} \right] \!+ \alpha_{K} f_{K}}{(1-\alpha_{K})\! \left[\sum_{k=1}^{K-1} \frac{\alpha_{k}}{1-\alpha_{K}} D_{k} \right] \!+ \alpha_{K} D_{K}} \!\stackrel{(a)}{\leq}\! \frac{f_{1}}{D_{1}}
\]
where (a) is from Fact $1$, noting that $f_{1}/D_{1} \geq f_{K}/D_{K}$ and
\[
\frac{f_{1}}{D_{1}} \geq \frac{\sum_{k=1}^{K-1} \frac{\alpha_{k}}{1-\alpha_{K}} f_{k}}{\sum_{k=1}^{K-1} \frac{\alpha_{k}}{1-\alpha_{K}} D_{k}},
\]
where the above holds by the induction assumption~\eqref{eq:405}.
\end{IEEEproof}

From Lemma~\ref{lem:401}, next we evaluate $f(\bQ(t))/\expect{T \mid \bQ(t)}$ for a given $\RRphi$ policy. Again we have $\expect{T\mid\bQ(t)} = \expect{T}$.

 In the special case $\bphi = \bm{0}$, we get $f(\bQ(t))/\expect{T \mid \bQ(t)} = 0$. Otherwise, fix some $\bphi\in\Phi$. For each active channel $n$ in $\bphi$,  we denote by $L_{n}^{\bphi}$ the amount of time the network stays with user $n$ in one round of $\RRphi$. It is shown in~\cite[Corollary $1$]{LaN10arXiv-channelmemory} that $L_{n}^{\bphi}$ has the probability distribution
\begin{equation} \label{eq:411}
L_{n}^{\bphi} = \begin{cases}
1 & \text{with prob. $1-\sP_{n, 01}^{(M(\bphi))}$} \\
j \geq 2 & \text{with prob. $\sP_{n, 01}^{(M(\bphi))} (\sP_{n, 11})^{(j-2)} \sP_{n, 10}$},
\end{cases}
\end{equation}
and
\begin{equation} \label{eq:409}
\expect{L_{n}^{\bphi}} = 1+ \frac{\sP_{n, 01}^{(M(\bphi))}}{\sP_{n, 10}}.
\end{equation}
It follows that under the $\RRphi$ policy we have
\begin{align*}
&\expect{T} = \expect{T \mid \bQ(t)} = \sum_{n: \phi_{n}=1} \expect{L_{n}^{\bphi}},  \\ 
&\expect{ \sum_{\tau=0}^{T-1} \mu_{n}(t+\tau) \mid \bQ(t)} = \begin{cases} 
	\expect{L_{n}^{\bphi}} - 1 & \text{if $\phi_{n}=1$} \\
	0 & \text{if $\phi_{n}=0$}
	\end{cases} 
\end{align*}
and thus
\begin{equation} \label{eq:408}
\frac{f(\bQ(t))}{\expect{T \mid \bQ(t)}} = \frac{\sum_{n=1}^{N} Q_{n}(t) \expect{L_{n}^{\bphi}-1} \phi_{n}}{\sum_{n=1}^{N} \expect{L_{n}^{\bphi}} \phi_{n}}.
\end{equation}

The above simplifications lead to the next network control algorithm that maximizes~\eqref{eq:401} in a frame-by-frame basis over all feasible admission and randomized round robin policies.

\underline{\bf Queue-dependent Round Robin for Network Utility}
\underline{\bf Maximization ($\qrrnum$):}

\begin{enumerate}
\item At the beginning of a transmission round, observe the current backlog vector $\bQ(t)$ and solve the convex program~\eqref{eq:209}-\eqref{eq:210}. Let $\br^{\qrr}(t) \triangleq (r_{n}^{\qrr}(t))_{n=1}^{N}$ be the optimal solution.

\item Let $\bphi^{\qrr}(t)$ be the maximizer of~\eqref{eq:408} over all  $\bphi\in\Phi$. If the resulting optimal objective is larger than zero, execute policy $\mathsf{RR}(\bphi^{\qrr}(t))$ for one round, with the channel ordering of least recently used first. Otherwise, idle the system for one slot. At the same time, admit $r_{n}^{\qrr}(t)$ packets to user $n$ in every slot of the current round. At the end of the round, go to Step 1).
\end{enumerate}

The most complex part of the $\qrrnum$ algorithm is to maximize~\eqref{eq:408} in Step 2, where in general all $(2^{N}-1)$  choices of vector $\bphi\in\Phi$ need to be examined, resulting in exponential complexity. In the special case that channels have the same transition probability matrix, the $\qrrnum$ algorithm reduces to a polynomial time policy, and the following  steps find the maximizer $\bphi^{\qrr}(t)$ of~\eqref{eq:408}:
\begin{enumerate}
\item Re-index $Q_{n}(t)$ so that $Q_{1}(t) \geq Q_{2}(t) \geq \ldots \geq Q_{N}(t)$.
\item For each $K\in\{1, \ldots, N\}$, compute
\begin{equation} \label{eq:410}
\frac{\sP_{01}^{(K)}}{K\left(\sP_{10}+\sP_{01}^{(K)}\right)} \sum_{k=1}^{K} Q_{k}(t)
\end{equation}
and let $K^{\qrr}$ be the maximizer of~\eqref{eq:410} over $K$.
\item If $\bQ(t)$ is the zero vector, let $\bphi^{\qrr}(t) = \bm{0}$. Otherwise, let $\bphi^{\qrr}(t)$ be the binary vector with the first $K^{\qrr}$ components being $1$ and $0$ otherwise.
\end{enumerate}

Another way to have an efficient $\qrrnum$ algorithm, especially when $N$ is large, is to restrict to a subset of $\RRphi$ policies. For example, consider those in every transmission round only serve $2$ or $0$ users. Although the associated new achievable region $\Lambdaint$ (can be found as a corollary of Theorem~\ref{thm:101}) will be smaller, the resulting $\qrrnum$ algorithm has polynomial time complexity because we only need to consider $N(N-1)/2$ choices of $\bphi$ in every round.


\section{Performance Analysis} \label{sec:analysis}

In the $\qrrnum$ policy, let $t_{k-1}$ and $T_{k}$ be the beginning and the duration of the $k$th transmission round. We have  $T_{k} = t_{k} - t_{k-1}$ and $t_{k} = \sum_{i=1}^{k} T_{i}$  for all $k\in\N$. Assume $t_{0} = 0$. Every $T_{k}$ is the length of a transmission round of some $\RRphi$ policy.  Define $\tmax$ as the length of a transmission round of the policy $\mathsf{RR}(\bm{1})$ that serves all channels in every round. Then for each $k\in\N$, we can show that $\tmax$ and $\tmax^{2}$ is stochastically larger than $T_{k}$ and $T_{k}^{2}$, respectively.\footnote{If $\RRphi$ for some $\bphi\in\Phi$ is used in $T_{k}$, the stochastic ordering between $\tmax$ and $T_{k}$ can be shown by noting that $T_{k} = \sum_{n: \phi_{n}=1}  L_{n}^{\bphi}$, where $L_{n}^{\bphi}$ is defined in~\eqref{eq:411}. Otherwise, we have $\bphi = \bm{0}$ and $T_{k}=1 \leq \tmax$.} As a result, we have
\begin{equation} \label{eq:711}
\expect{T_{k}}\leq \expect{\tmax} < \infty , \quad \expect{T_{k}^{2}} \leq \expect{\tmax^{2}} <\infty.
\end{equation}
The next theorem shows the performance of the $\qrrnum$ algorithm.

\begin{thm} \label{thm:401}
Let $\by(t) = (y_{n}(t))_{n=1}^{N}$ be the vector of served packets for each user in slot $t$; $y_{n}(t) = \min[Q_{n}(t), \mu_{n}(t)]$.
Define constant $B \triangleq N\expect{\tmax^{2}}$. Then for any given positive control parameter $V_{g}>0$, the $\qrrnum$ algorithm stabilizes the network and yields average network utility satisfying
\begin{equation} \label{eq:701}
\liminf_{t\to\infty} g\left(\frac{1}{t} \sum_{\tau=0}^{t-1} \expect{\by(\tau)}\right) \geq g(\boy^{*}) - \frac{B}{V_{g}},
\end{equation}
where $g(\boy^{*})$ is the optimal network utility and the solution to the constrained restless bandit problem~\eqref{eq:104}-\eqref{eq:105}. By taking $V_{g}$ sufficiently large, the $\qrrnum$ algorithm achieves network utility arbitrarily close to the optimal $g(\boy^{*})$, and thus solves~\eqref{eq:104}-\eqref{eq:105}.
\end{thm}

\begin{IEEEproof}[Proof of Theorem~\ref{thm:401}]
Analyzing the performance of the $\qrrnum$ algorithm relies on comparing it to a near-optimal feasible solution. We will adopt the approach in~\cite{NML08} but generalize it to a frame-based analysis.

For some $\epsilon >0$, consider the \emph{$\epsilon$-constrained} version of~\eqref{eq:104}-\eqref{eq:105}:
\begin{align}
\text{maximize:} \quad & g(\boy) \label{eq:301} \\
\text{subject to:} \quad &  \boy \in\Lambdaint(\epsilon) \label{eq:302}
\end{align}
where $\Lambdaint(\epsilon)$ is the achievable region $\Lambdaint$ stripping an ``$\epsilon$-layer''  off the boundary:
\[
\Lambdaint(\epsilon) \triangleq \{ \boy \mid \boy + \epsilon \bm{1} \in \Lambdaint\},
\]
where $\bm{1}$ is an all-one vector. Notice that $\Lambdaint(\epsilon) \to \Lambdaint$ as $\epsilon\to 0$. Let $\boy^{*}(\epsilon) = (\oy_{n}^{*}(\epsilon))_{n=1}^{N}$ and $\boy^{*} = (\oy_{n}^{*})_{n=1}^{N}$ be the optimal solution to the $\epsilon$-constrained problem~\eqref{eq:301}-\eqref{eq:302} and the constrained restless bandit problem~\eqref{eq:104}-\eqref{eq:105}, respectively. For simplicity, we assume $\boy^{*}_{\epsilon} \to \boy^{*}$ as $\epsilon \to 0$.\footnote{This property is proved in a similar case in~\cite[Ch. $5.5.2$]{Nee03thesis}.}

From Corollary~\ref{cor:301}, there exists a randomized round robin that yields the throughput vector $\boy^{*}_{\epsilon} + \epsilon \bm{1}$ (note that $\boy^{*}_{\epsilon} + \epsilon \bm{1}\in \Lambdaint$), and we denote this policy by $\randrr^{*}_{\epsilon}$. Let $T^{*}_{\epsilon}$ denotes the length of one transmission round under $\randrr^{*}_{\epsilon}$. Then we have for each $n\in\{1, \ldots, N\}$
\begin{equation} \label{eq:304}
\expect{\sum_{\tau=0}^{T^{*}_{\epsilon}-1} \mu_{n}(t+\tau) \mid \bQ(t)} \geq (\oy_{n}^{*}(\epsilon) +\epsilon) \expect{T^{*}_{\epsilon}}
\end{equation}
from renewal reward theory. That is, we may consider a renewal reward process where renewal epochs are time instants at which $\randrr^{*}_{\epsilon}$ starts a new round of transmission (with renewal period $T^{*}_{\epsilon}$), and rewards are the allocated service rates $\mu_{n}(t+\tau)$.\footnote{We note that this renewal reward process is defined solely with respect to the service policy $\randrr_{\epsilon}^{*}$, and the network state needs not renew itself at the renewal epochs.} Then the average service rate is simply the average sum reward over a renewal period divided by the average renewal duration $\expect{T^{*}_{\epsilon}}$. Then~\eqref{eq:304} holds because this average service rate is greater than or equal to $(\oy^{*}(\epsilon)+\epsilon)$.

Combining $\randrr_{\epsilon}^{*}$ with the admission policy $\sigma^{*}$ that sets $r_{n}(t+\tau) = \oy_{n}^{*}(\epsilon)$ for all $n$ and $\tau\in\{0, \ldots, T^{*}_{\epsilon}-1\}$,\footnote{Since the throughput vector $\boy^{*}(\epsilon) = (\oy^{*}_{n}(\epsilon))_{n=1}^{N}$ is achievable in $\Lambdaint(\epsilon)$, each component $\oy_{n}^{*}(\epsilon)$ must be less than or equal to the stationary probability $\pi_{n, \sON} \leq 1$, and thus is a feasible choice of $r_{n}(t)$.} we get
\begin{align}
&f_{\epsilon}^{*}(\bQ(t)) \geq \expect{T^{*}_{\epsilon}} \sum_{n=1}^{N} Q_{n}(t) (\oy_{n}^{*}(\epsilon) +\epsilon)  \label{eq:310} \\
&g_{\epsilon}^{*}(\bQ(t))  = \expect{T^{*}_{\epsilon}} \left[V_{g}\,g(\boy^{*}(\epsilon)) - \sum_{n=1}^{N} Q_{n}(t) \oy_{n}^{*}(\epsilon) \right]  \label{eq:311}
\end{align}
where~\eqref{eq:310}\eqref{eq:311} are $f(\bQ(t))$ and $g(\bQ(t))$ evaluated under $\randrr_{\epsilon}^{*}$ and $\sigma^{*}$, respectively.

Since the $\qrrnum$ policy maximizes~\eqref{eq:401}, evaluating~\eqref{eq:401} under both $\qrrnum$ and the policy $(\randrr_{\epsilon}^{*}, \sigma^{*})$ yields
\begin{equation} \label{eq:312}
\begin{split}
&f_{\qrrnum}(\bQ(t_{k})) + g_{\qrrnum}(\bQ(t_{k})) \\
& \geq \expect{T_{k+1}\mid \bQ(t_{k})} \frac{f_{\epsilon}^{*}(\bQ(t_{k})) + g_{\epsilon}^{*}(\bQ(t_{k}))}{\expect{T_{\epsilon}^{*}}} \\
& \stackrel{(a)}{\geq} \expect{T_{k+1} \mid \bQ(t_{k})} \left[V_{g}\,g(\boy^{*}(\epsilon)) + \epsilon \sum_{n=1}^{N} Q_{n}(t_{k}) \right] \\
& = \expect{
T_{k+1} \left(V_{g}\,g(\boy^{*}(\epsilon)) + \epsilon \sum_{n=1}^{N} Q_{n}(t_{k}) \right) \mid \bQ(t_{k})
},
\end{split}
\end{equation}
where (a) is from~\eqref{eq:310}\eqref{eq:311}. The drift-minus-utility bound~\eqref{eq:309} under the $\qrrnum$ policy in the $(k+1)$th round of transmission yields
\begin{equation} \label{eq:406}
\begin{split}
&\Delta_{T_{k+1}}(\bQ(t_{k})) - V_{g} \expect{ \sum_{\tau=0}^{T_{k+1}-1} g(\br(t_{k}+\tau)) \mid \bQ(t_{k})} \\
&\quad \leq B - f_{\qrrnum}(\bQ(t_{k})) - g_{\qrrnum}(\bQ(t_{k})) \\
&\quad \stackrel{(a)}{\leq} B - \expect{
T_{k+1} \left(V_{g}\,g(\boy^{*}(\epsilon)) + \epsilon \sum_{n=1}^{N} Q_{n}(t_{k}) \right) \mid \bQ(t_{k})
}
\end{split}
\end{equation}
where (a) is from~\eqref{eq:312}. Taking expectation over $\bQ(t_{k})$ in~\eqref{eq:406} and summing it over $k\in\{0, \ldots, K-1\}$, we get
\begin{equation} \label{eq:407}
\begin{split}
&\expect{L(\bQ(t_{K}))} - \expect{L(\bQ(t_{0}))} - V_{g} \expect{ \sum_{\tau=0}^{t_{K}-1} g(\br(\tau))} \\
&\leq BK - V_{g} g(\boy^{*}(\epsilon)) \expect{t_{K}} - \epsilon\,\expect{ \sum_{k=0}^{K-1} T_{k+1} \sum_{n=1}^{N} Q_{n}(t_{k})}.
\end{split}
\end{equation}
Since $Q_{n}(\cdot)$ and $L(\bQ(\cdot))$ are nonnegative and $\bQ(t_{0})=\bm{0}$, ignoring all backlog-related terms in~\eqref{eq:407} yields
\begin{equation} \label{eq:414}
\begin{split}
- V_{g} \expect{ \sum_{\tau=0}^{t_{K}-1} g(\br(\tau))} &\leq BK - V_{g} g(\boy^{*}(\epsilon)) \expect{t_{K}} \\
&\stackrel{(a)}{\leq} B\expect{t_{K}} - V_{g} g(\boy^{*}(\epsilon)) \expect{t_{K}}
\end{split}
\end{equation}
where (a) uses $t_{K} = \sum_{k=1}^{K} T_{k} \geq K$. Dividing~\eqref{eq:414} by $V_{g}$ and rearranging terms, we get
\begin{equation} \label{eq:415}
\expect{ \sum_{\tau=0}^{t_{K}-1} g(\br(\tau))} \geq \left( g(\boy^{*}(\epsilon)) - \frac{B}{V_{g}} \right) \expect{t_{K}}.
\end{equation}
Recall from Section~\ref{sec:drift} that $B$ is an unspecified constant satisfying $B\geq N \expect{T_{k}^{2} \mid \bQ(t)}$. From~\eqref{eq:711} it suffices to define
$
B \triangleq N \expect{\tmax^{2}}
$.

In $\qrrnum$, let $K(t)$ denote the number of transmission rounds ending before time $t$. Using $t_{K(t)} \leq t < t_{K(t)+1}$, we have $0 \leq t - \expect{t_{K(t)}} \leq \expect{t_{K(t)+1}-t_{K(t)}} = \expect{T_{K(t)+1}}$. Dividing the above by $t$ and passing $t\to\infty$, we get
\begin{equation} \label{eq:418}
\lim_{t\to\infty} \frac{t - \expect{t_{K(t)}}}{t} = 0.
\end{equation}
Next, the expected sum utility over the first $t$ slots satisfies
\begin{equation} \label{eq:416}
\begin{split}
&\sum_{\tau=0}^{t-1} \expect{g(\br(\tau))} = \expect{ \sum_{\tau=0}^{t_{K(t)}-1} g(\br(\tau))} + \expect{\sum_{\tau=t_{K(t)}}^{t-1} g(\br(\tau)) } \\ 
&\stackrel{(a)}{\geq} \left[ g(\boy^{*}(\epsilon)) - \frac{B}{V_{g}} \right] \expect{t_{K(t)}}  \\
&= \left[ g(\boy^{*}(\epsilon)) - \frac{B}{V_{g}} \right] t - \left[g(\boy^{*}(\epsilon)) - \frac{B}{V_{g}} \right] \left(t- \expect{t_{K(t)}}\right),
\end{split}
\end{equation}
where (a) uses~\eqref{eq:415} and that $g(\cdot)$ is nonnegative. Dividing~\eqref{eq:416} by $t$, taking a $\liminf$ as $t\to\infty$ and using~\eqref{eq:418}, we get
\begin{equation} \label{eq:417}
\liminf_{t\to\infty} \frac{1}{t} \sum_{\tau=0}^{t-1} \expect{g(\br(\tau))} \geq g(\boy^{*}(\epsilon)) - \frac{B}{V_{g}}.
\end{equation}
Using Jensen's inequality and the concavity of $g(\cdot)$, we get
\begin{equation} \label{eq:420}
\liminf_{t\to\infty} \frac{1}{t} \sum_{\tau=0}^{t-1} \expect{g(\br(\tau))} \leq \liminf_{t\to\infty}  g\left( \bor^{(t)} \right),
\end{equation}
where we define  the average admission data vector:
\begin{equation} \label{eq:712}
\bor^{(t)} \triangleq (\oor_{n}^{(t)})_{n=1}^{N}, \quad \oor_{n}^{(t)} \triangleq \frac{1}{t} \sum_{\tau=0}^{t-1} \expect{\br(\tau)}. 
\end{equation}
Combining~\eqref{eq:417}\eqref{eq:420} yields
\[
\liminf_{t\to\infty}  g\left( \bor^{(t)} \right) \geq g(\boy^{*}(\epsilon)) - \frac{B}{V_{g}},
\]
which holds for any sufficiently small $\epsilon$. Passing $\epsilon\to 0$ yields
\begin{equation} \label{eq:412}
\liminf_{t\to\infty}  g\left(  \bor^{(t)}  \right) \geq g(\boy^{*}) - \frac{B}{V_{g}}.
\end{equation}

Finally, we show the network is stable, and as a result
\begin{equation} \label{eq:702}
\liminf_{t\to\infty}  g\left( \boy^{(t)}  \right) \geq \liminf_{t\to\infty}  g\left(  \bor^{(t)}  \right),
\end{equation}
where $\boy^{(t)} = (\oy_{n}^{(t)})_{n=1}^{N}$ is defined similarly as $\bor^{(t)}$. Then combining~\eqref{eq:412}\eqref{eq:702} finishes the proof. To prove stability, ignoring the first, second, and fifth term in~\eqref{eq:407} yields
\begin{equation} \label{eq:419}
\begin{split}
\epsilon \expect{ \sum_{k=0}^{K-1} T_{k+1} \sum_{n=1}^{N} Q_{n}(t_{k})} &\leq BK + V_{g} \expect{\sum_{\tau=0}^{t_{K}-1} g(\br(\tau))} \\
& \stackrel{(a)}{\leq} K\left(B + V_{g} \gmax \expect{\tmax}\right)
\end{split}
\end{equation}
where we define $\gmax \triangleq g(\bm{1}) < \infty$ as the maximum value of $g(\cdot)$ (since $g(\cdot)$ is nondecreasing), and (a) uses
\[
g(\br(\tau)) \leq \gmax, \quad \expect{t_{K}} = \sum_{k=1}^{K} \expect{T_{k}} \leq K \expect{T_{\text{max}}}.
\]

Dividing~\eqref{eq:419} by $K\epsilon$, taking a $\limsup$ as $K\to\infty$, and using $T_{k+1} \geq 1$, we get
\begin{equation} \label{eq:421}
\limsup_{K\to\infty} \frac{1}{K} \expect{ \sum_{k=0}^{K-1} \sum_{n=1}^{N} Q_{n}(t_{k})} \leq \frac{B+ V_{g} \gmax \expect{T_{\text{max}}}}{\epsilon} < \infty.
\end{equation}
Equation~\eqref{eq:421} shows that the average backlog is bounded when sampled at time instants $\{t_{k}\}$. This property is enough to conclude that the average backlog over the whole time horizon is bounded, namely~\eqref{eq:706} holds and the network is stable. It is because the length of each transmission round $T_{k}$ has a finite second moment and the maximum amount of data admitted to each user in every slot is at most $1$; see~\cite[Lemma~$13$]{LaN10arXiv-channelmemory} for a detailed proof.

It remains to show network stability leads to~\eqref{eq:702}. Recall that $y_{n}(\tau) = \min[Q_{n}(\tau), \mu_{n}(\tau)]$ is the number of user-$n$ packets served in slot $\tau$, and~\eqref{eq:202} is equivalent to
\begin{equation} \label{eq:703}
Q_{n}(\tau+1) = Q_{n}(\tau) - y_{n}(\tau) + r_{n}(\tau).
\end{equation}
Summing~\eqref{eq:703} over $\tau\in\{0, \ldots, t-1\}$, taking an expectation and dividing it by $t$, we get
\begin{equation} \label{eq:704}
\frac{\expect{Q_{n}(t)}}{t} = \oor_{n}^{(t)} - \oy_{n}^{(t)},
\end{equation}
where $\oor_{n}^{(t)}$ is defined in~\eqref{eq:712} and $\oy_{n}^{(t)}$ is defined similarly.
From~\cite[Theorem $4(c)$]{Nee10arXiv-stability}, the stability of $Q_{n}(t)$ and~\eqref{eq:704} result in that for each $n$:
\begin{equation} \label{eq:707}
\limsup_{t\to\infty} \frac{\expect{Q_{n}(t)}}{t} = \limsup_{t\to\infty} \left( \oor_{n}^{(t)} - \oy_{n}^{(t)} \right) = 0.
\end{equation}
Noting that $g(\cdot)$ is bounded, there exists a convergent subsequence of $g(\boy^{(t)})$ indexed by $\{t_{i}\}_{i=1}^{\infty}$ such that
\begin{equation} \label{eq:708}
\lim_{i\to\infty} g\left(\boy^{(t_{i})}\right) =  \liminf_{t\to\infty} g\left(\boy^{(t)}\right).
\end{equation}
By iteratively finding a convergent subsequence of $\{\oor_{n}^{(t_{i})}\}_{i=1}^{\infty}$ for each $n$ (noting that $\oor_{n}^{(t)}$ is bounded for all $n$ and $t$), there exists a subsequence $\{t_{k}\} \subset \{t_{i}\}$ such that $\{\bor^{(t_{k})}\}_{k=1}^{\infty}$ converges as $k\to\infty$. From~\eqref{eq:707} and that $\limsup \{z_{n}\}$ is the supremum of all limit points of a sequence $\{z_{n}\}$, we get
\begin{equation} \label{eq:709}
\lim_{k\to\infty} (\oor_{n}^{(t_{k})} - \oy_{n}^{(t_{k})}) \leq 0 \Rightarrow \lim_{k\to\infty} \oor_{n}^{(t_{k})} \leq \lim_{k\to\infty} \oy_{n}^{(t_{k})}, \quad \forall n.
\end{equation}
It follows that
\begin{align*}
\liminf_{t\to\infty} g\left(\boy^{(t)}\right) &\stackrel{(a)}{=} \lim_{k\to\infty} g\left(\boy^{(t_{k})}\right) \\ &\stackrel{(b)}{\geq} \lim_{k\to\infty} g\left(\bor^{(t_{k})}\right)  \stackrel{(c)}{\geq} \liminf_{t\to\infty} g\left(\bor^{(t)}\right),
\end{align*}
where (a) is from~\eqref{eq:708}, (b) uses~\eqref{eq:709} and that  $g(\cdot)$ is continuous and nondecreasing, and (c) uses that $\liminf \{z_{n}\}$ is the infimum of limit points of $\{z_{n}\}$.
\end{IEEEproof}

\section{Conclusions}

We have provided a theoretical framework to do network utility maximization over partially observable Markov $\sON$/$\sOFF$ channels. 
The performance and control decisions in such networks are constrained by the limiting channel probing capability and delayed/uncertain channel state information, but can be improved by taking advantage of channel memory. Overall, to attack such problems we need to solve (at least approximately) high-dimensional restless bandit problems with a general functional objective, which are difficult to analyze using existing tools such as Whittle's index theory or Markov decision theory. In this paper we propose a new methodology to solve such problems by combining an achievable region approach from mathematical programming and the powerful Lyapunov optimization theory. The key idea is to first identify a good constrained performance region rendered by stationary policies, and then solve the problem only over the constrained region, serving as an approximation to the original problem. While a constrained performance region is constructed in~\cite{LaN10arXiv-channelmemory}, in this paper using a novel frame-based variable-length Lyapunov drift argument, we can solve the original problem over the constrained region by constructing queue-dependent greedy algorithms that stabilize the network with near-optimal utility. It will be interesting to see how the Lyapunov optimization theory can be extended and used to attack other sequential decision making problems as well as stochastic network optimization problems with limited channel probing and delayed/uncertain channel state information.

\bibliographystyle{IEEEtran}
\bibliography{/Users/chihping/Desktop/bibliography/IEEEabrv,/Users/chihping/Desktop/bibliography/myabrv,/Users/chihping/Desktop/bibliography/mypaperbib}
\end{document}